\documentclass{article}     
\usepackage{latexsym,amsthm,amssymb,amsmath}
\usepackage{enumerate}
\usepackage{proof,logic}

\newtheorem{theorem}{Theorem}

\title{On non-eliminability of the cut rule and the roles of associativity and distributivity in non-commutative substructural logics}
\author{
Takeshi Ueno\thanks{Rakuno-Gakuen University, Japan, t-ueno@rakuno.ac.jp}
\and
Koji nakaogawa\thanks{Department of Philosophy, Hokkaido University, Japan, koji@logic.let.hokudai.ac.jp}
\and
Osamu Watari\thanks{Hokkaido Automotive Engineering College, Japan, watari@haec.ac.jp}
}

\date{}

\begin{document}
\maketitle

\begin{abstract}
We introduce a sequent calculus \FL{}{}' ,which has at most one formula on the right side of sequent, and which excludes three structural inference rules, i.e. contraction, weakening and exchange. Our formulations of the inference rules of  \FL{}{}' are based on the results and considerations carried out in our previous papers on how to formulate Gentzen-style natural deduction for non-commutative substructural logics.

Our present formulation  \FL{}{}' of sequent system for non-commutative substructural logic, which has no structural rules, has the same proof strength as the ordinary and standard sequent calculus  \FL{}{} (Full Lambek), which is often called Full Lambek calculus, i.e., the basic sequent calculus for all other substructural logics. For the standard \FL{}{} (Full Lambek), we use Ono's formulation.

Although our  \FL{}{}' and the standard formulation  \FL{}{} (Full Lambek) are equivalent, there is a subtle difference in the left rule of implication. In the standard formulation, two parameters $\Gamma_1$ and $\Gamma_2$(resp.), each of which is just an finite sequence of arbitrary formulas, appear on the left and right side (resp.) of a formula appearing on the left side of the sequent on the upper left side the left rule $\imply$\;(which corresponds to $\imply'$ in \FL{}{}') . On the other hand, there is no such parameter on the left side of the sequent on the upper left side in the left rule for $\imply'$ of our system \FL{}{}. In our system  \FL{}{}', $\Gamma_1$ is always empty, and only $\Gamma_2$ is allowed to occur in the left rule for $\imply'$ (similar differences occur in the multiplicative conjunction, additive conjunction and additive disjunction). This subtle difference between our system  \FL{}{}' and the standard system  \FL{}{} (Full Lambek) matters deeply, for we are led to a construction of proof-figures in  \FL{}{}', which show how the associativity of multiplicative conjunction and the distributivity of multiplicative conjunction over additive disjunction are involved in the eliminations of the cut rule in those proofs. We clarify and specify how associativity and distributivity are related to the non-eliminability of an application of the cut rule in those proof-figures of  \FL{}{}'.

\end{abstract}

%
%
%
%
\section{Introduction}


The situation surrounding the syntactic aspects of non-commutative substructural logics does not seem to be fully clarified. In particular, the process of eliminating applications of the cut rule in a given proof-figure of intuitionistic sequent system for $\FL'{}{}$ (defined below) , where \FL{}{} stands for "Full Lambek" and is the most basic system for substructural logics, sometimes succeeds and terminates, and some other times, does not succeed and does not terminate to produce a proof-figure which contains no applications of the cut rule. Indeed, it depends on the subtlety of where one is allowed to place parameters (side formulas) in some of the inference rules of $\FL'{}{}$ . 
     
In the present paper, we introduce a system of inference rules of intuitionistic ("left heavy") sequent calculus for substructural logic, \FL{}{}', which lacks all structural inference rules, namely, exchange, weakening, and contraction rules. Furthermore, the parameters in its inference rules are placed in such a way that their positions reflect the "natural order" of non-cancelled hypothesis in the (Gentzen-style) natural deduction for non-commutative substructural logic. Using our system $\FL'{}{}$ , we will show how the associative law for multiplicative conjunction and the distributive law of multiplicative conjunction over additive disjunction are entangled in the elimination process of applications of the cut rule. Analysis the relevancy of these two rules as to the cut elimination process has become possible to us, for we fixed the positions of parameters in the inference rules of $\FL'{}{}$  according to our analysis of normalization procedures in Gentzen-style natural deduction for non-commutative substructural logic. (This paper does not assume the knowledge of our previous papers on Gentzen-style natural deductions for substructural logics.)  


%
%
\section{Language  ${\cal L}$ and its Formulas}

Our language  ${\cal L}$ has propositional constant symbols $A, B, C, \cdots$ .
As for logical connectives, it has
the implication symbols $\imply $, $\coimply $,
the negation symbols $\lnot{}$, $\colnot{}$,
the multiplicative conjunction symbol $\tprod$,
the additive conjunction symbol $\dprod$,
the additive disjunction symbol $\dsum$.
In  ${\cal L}$, there are constant symbols $\tunit$ to denote
the unit element for the multiplicative conjunction,
$\tzero$ to denote the unit element for the multiplicative
disjunction which is not introduced in this paper,
$\dunit$ to denote
the unit element for the additive conjunction, and
$\dzero$ to denote the unit element for the additive disjunction.

The formulas of  ${\cal L}$ are defined inductively as a finite sequence
of these symbols together with parenthesizes.

%
%
%
\section{Sequent Calculus \FL{}{}'\;\;(our formulation)}

The sequent of the language ${\cal L}$ have the following form

\begin{center}
$\Gamma$ $\to$ $\Delta$ .
\end{center}
The left hand side of a sequent may be empty. The
right hand side of a sequent is either empty or
consists of a single formula.
To specify the element of \ $\Gamma$ \ and \
$\Delta$ \, we write
$$
\gamma_{0}, \cdots, \gamma_{m-1} \to\delta.
$$
When both sides of a sequent are empty, we write
$$
\to
$$

Next, we introduce  the sequent calculus
 $\FL{}{}'$ as follows.
We say that \ $(\Gamma_1 \to X_1, \Gamma_2  \to X_2,
\cdots , \Gamma_n \to X_n \ / \ \Gamma \to X)$ \
is an instance of a certain inference rule if it
has the form indicated by the corresponding figure.
If \ $(\Gamma_1 \to X_1, \Gamma_2  \to X_2, \cdots ,
\Gamma_n \to X_n \ / \ \Gamma \to X)$ \ is an instance
of an inference rule $\alpha$ , we call \ $\Gamma_i
\to X_i$ \ the $i$-th upper sequent of $\alpha$ , and \
$\Gamma \to X$ \ the lower sequent of \ $\alpha$.
(The origin of $\FL{}{}$ goes back to a classical paper
written by J. Lambek in 1950's. 
Our presentation of $\FL{}{}'$ is based on Ono[4].), but ours is different from his in important places.


In our LF', the positions of parameters in its inference rules are determined according to the un-cancelled hypothesis of Gentzen-style natural deduction for non-commutative substructural logic. The reader should take note that the positions of parameters in inferences rules of our \FL{}{}' are different from those of Ono's. 


\begin{itemize}
\item Axioms and rule for logical constants:
$$
\arraycolsep1em\def\arraystretch{2}
\begin{array}{cccc}
A  \drv A & \\
   \drv \tunit & \tzero \drv   \\
\G \drv \dunit & \dzero, \G \drv C \\
\infer[\tunit \Weakening]{\tunit, \G \drv C}{\G \drv C} &
\infer[\tzero \Weakening]{\G \drv \tzero}{\G \drv }
\end{array}
$$

\item Structural inference rule:

$$
\infer[\Cut]{\G_{2}, \G_{1}, \G_{3} \drv C}{
  \G_{1} \drv A & \G_{2}, A, \G_{3} \drv C}
$$

\item Logical inference rule:

$$
\arraycolsep1em\def\arraystretch{2}
\begin{array}{cc}
\infer[\LR{\imply}]{\G_{1}, A \imply B, \G_{2} \drv C}{
  \G_{1} \drv A & B, \G_{2} \drv C}
&
\infer[\RR{\imply}]{\G \drv A \imply B}{
  A, \G \drv B}
\\
\infer[\LR{\coimply}]{A \coimply B, \G_{1}, \G_{2} \drv C}{
  \G_{1} \drv A & B, \G_{2} \drv C}
&
\infer[\RR{\coimply}]{\G \drv A \coimply B}{
  \G, A \drv B}
\\
\infer[\LR{\lnot{}}]{\G, \lnot{A} \drv}{
  \G \drv A}
&
\infer[\RR{\lnot{}}]{\G \drv \lnot{A}}{
  A, \G \drv }
\\
\infer[\LR{\colnot{}}]{\colnot{A}, \G \drv}{
  \G \drv A}
&
\infer[\RR{\colnot{}}]{\G \drv \colnot{A}}{
  \G,A \drv }
\\
\infer[\LR{\tprod}]{A \tprod B, \G \drv C}{
  A, B, \G \drv C}
&
\infer[\RR{\tprod}]{\G_{1}, \G_{2} \drv A \tprod B}{
  \G_{1} \drv A & \G_{2} \drv B}
\\
\begin{array}{c}
\infer[\LR{\dprod_{1}}]{A \dprod B, \G \drv C}{
  A, \G \drv C}
\\
\infer[\LR{\dprod_{2}}]{A \dprod B, \G \drv C}{
  B, \G \drv C}
\end{array}
&
\infer[\RR{\dprod}]{\G \drv A \dprod B}{
  \G \drv A & \G \drv B}
\\
\infer[\LR{\dsum}]{A \dsum B, \G \drv C}{
  A, \G \drv C & B, \G \drv C}
&
\begin{array}{c}
\infer[\RR{\dsum_{1}}]{\G \drv A \dsum B}{
  \G \drv A}
\\
\infer[\RR{\dsum_{2}}]{\G \drv A \dsum B}{
  \G \drv B}
\end{array}
\end{array}
$$

\end{itemize}

%
%
%
\section{Sequent Calculus \FL{}{}\;\;(Ono's formulation)}

The reader should be warned that the $\imply$ of \FL{}{}' corresponds with $\imply'$ of Ono's, and
$\imply'$ of \FL{}{}' corresponds with $\imply$ of Ono's. 
$\lnot$ of \FL{}{}' corresponds with $\lnot'$ of Ono's, and
$\lnot'$ of \FL{}{}' corresponds with $\lnot$ of Ono's.

\begin{itemize}
\item Axioms and rules for logical constants:
$$
\arraycolsep1em\def\arraystretch{2}
\begin{array}{cccc}
A  \drv A & \\
   \drv \tunit & \tzero \drv   \\
\G \drv \dunit & \G_{1}, \dzero, \G_{2} \drv C \\
\infer[\tunit \Weakening]{\G_{1}, \tunit, \G_{2} \drv C}{\G_{1}, \G_{2} \drv C} &
\infer[\tzero \Weakening]{\G \drv \tzero}{\G \drv }
\end{array}
$$

\item Structural inference rules:

$$
\infer[\Cut]{\G_{2}, \G_{1}, \G_{3} \drv C}{
  \G_{1} \drv A & \G_{2}, A, \G_{3} \drv C}
$$

\item Logical inference rules:

$$
\arraycolsep1em\def\arraystretch{2}
\begin{array}{cc}
\infer[\LR{\coimply}]{\G_{2}, \G_{1}, A \coimply B, \G_{3} \drv C}{
  \G_{1} \drv A & \G_{2}, B, \G_{3} \drv C}
&
\infer[\RR{\coimply}]{\G \drv A \coimply B}{
  A, \G \drv B}
\\
\infer[\LR{\imply}]{\G_{2}, A \imply B, \G_{1}, \G_{3} \drv C}{
  \G_{1} \drv A & \G_{2}, B, \G_{3} \drv C}
&
\infer[\RR{\imply}]{\G \drv A \imply B}{
  \G, A \drv B}
\\
\infer[\LR{\colnot{}}]{\G, \colnot{A} \drv}{
  \G \drv A}
&
\infer[\RR{\colnot{}}]{\G \drv \colnot{A}}{
  A, \G \drv }
\\
\infer[\LR{\lnot{}}]{\lnot{A}, \G \drv}{
  \G \drv A}
&
\infer[\RR{\lnot{}}]{\G \drv \lnot{A}}{
  \G,A \drv }
\\
\infer[\LR{\tprod}]{\G_{1}, A \tprod B, \G_{2} \drv C}{
  \G_{1}, A, B, \G_{2} \drv C}
&
\infer[\RR{\tprod}]{\G_{1}, \G_{2} \drv A \tprod B}{
  \G_{1} \drv A & \G_{2} \drv B}
\\
\begin{array}{c}
\infer[\LR{\dprod_{1}}]{\G_{1}A \dprod B, \G_{2} \drv C}{
  \G_{1}, A, \G_{2} \drv C}
\\
\infer[\LR{\dprod_{2}}]{\G_{1}, A \dprod B, \G_{2},  \drv C}{
  \G_{1}, B, \G_{2} \drv C}
\end{array}
&
\infer[\RR{\dprod}]{\G \drv A \dprod B}{
  \G \drv A & \G \drv B}
\\
\infer[\LR{\dsum}]{\G_{1}, A \dsum B, \G_{2} \drv C}{
  \G_{1}, A, \G_{2} \drv C & \G_{1}, B, \G_{2}, \drv C}
&
\begin{array}{c}
\infer[\RR{\dsum_{1}}]{\G \drv A \dsum B}{
  \G \drv A}
\\
\infer[\RR{\dsum_{2}}]{\G \drv A \dsum B}{
  \G \drv B}
\end{array}
\end{array}
$$

\end{itemize}

%
%
%

\section{Equivalence of \FL{}{} and  $\FL{}{}'$ }

%
%
\begin{theorem}[Equivalence of $\FL{}{}$ and \FL{}{}']
\mbox{} \\
Let \ $\phi$ \ be a formula of the language  ${\cal L}$.
Let \ $\G$ \ be a list of formulas of  ${\cal L}$.
Then, \ the sequent
\ $\G \drv \phi$ \ is provable in  $\FL{}{}$
\ if and only if the sequent
\ $\G \drv \phi$ \ is provable in  $\FL{}{}'$.
\end{theorem}

\paragraph{{\it Proof}}

First, we prove  that if sequent $\G \drv \phi$
is provable with proof $\Pi$ in $\FL{}{}$,
then sequent $\G \drv \phi$ is provable with a proof $\S$ in \FL{}{}'.
To prove this direction, we use induction on the number $\sharp (\Pi)$ of
the applications of inference rules in the proof $\Pi$.

If $\sharp (\Pi)$ is zero, $\G \drv \phi$ must be
an axiom of $\FL{}{}$.

The axiom 
$ \G \drv \dunit$ ,\;  $\G_{1}, \dzero, \G_{2} \drv C $
of  $\FL{}{}$ is provable in  $\FL{}{}'$.
This is shown by the following proof figure of  $\FL{}{}'$.

$$
\infer[\Cut]{\alpha_1, \alpha_2, \dzero , \G_2 \drv C}{
  \dzero , \G_2 \drv \alpha_1\tprod\alpha_2\imply C
  &
  \infer[\LR{\imply}]{\alpha_1, \alpha_2, \alpha_1\tprod\alpha_2\imply C}{
    \infer[\RR{\tprod}]{\alpha_1, \alpha_2 \drv \alpha_1\tprod\alpha_2}{
      \alpha_1 \drv \alpha_1
      &
      \alpha_2 \drv \alpha_2
    }
    &
    C\drv C
  }
}
$$


The following proof figure shows that the axiom (rule)  \tunit \Weakening \; of $\FL{}{}$ 
is provable in  $\FL{}{}'$.

$$
\infer[\Cut]{\alpha_1, \alpha_2, \tunit , \G_2 \drv C}{
  \infer[\tunit\Weakening]{\tunit, \G_2 \drv \alpha_1\tprod\alpha_2\imply C}{
    \infer[\RR{\imply}]{\G_2 \drv \alpha_1\tprod\alpha_2\imply C}{
      \infer[\LR{\imply}]{\alpha_1\tprod\alpha_2, \G_2 \drv C}{
      \alpha_1, \alpha_2, \G_2 \drv C
      }
    }
  }
  &
  \infer[\LR{\imply}]{\alpha_1, \alpha_2, \alpha_1\tprod\alpha_2\imply C}{
    \infer[\RR{\tprod}]{\alpha_1, \alpha_2 \drv \alpha_1\tprod\alpha_2}{
      \alpha_1 \drv \alpha_1
      &
      \alpha_2 \drv \alpha_2
    }
    &
    C\drv C
  }
}
$$

Now, we assume the theorem for $\sharp (\Pi) < n$,
and prove it for $\sharp \ (\Pi) = n$.

Our proof is divided into cases, depending on which inference rule
is used as the ``bottom'' inference rule in $\Pi$.
%
%

Without loss of generality, we assume that $\G_1$ consits of just  $\alpha_1$,\;and \;$\alpha_2$.

%
\vspace*{2em}
\begin{description}

\item[Case 1]%
The bottom inference rule in $\Pi$ is$\LR{\tprod}$ \\[0.5em]
i.e.\;\;
$
\infer[\LR{\tprod}]{\alpha_1, \alpha_2, A \tprod B, \G_{2} \drv C}{
  \alpha_1, \alpha_2, A, B, \G_{2} \drv C}
$ \\[1em]
Then we can construct a proof of \;$\alpha_1, \alpha_2, A \tprod B, \G_{2} \drv C$ in $\FL{}{}'$ as follows:
$$
\infer[\Cut]{\alpha_1, \alpha_2, A\tprod B, \G_2 \drv C}{
  \infer[\LR{\tprod}]{A\tprod B, \G_2 \drv \alpha_1\tprod\alpha_2\imply C}{
    \infer[\RR{\imply}]{A, B, \G_2 \drv \alpha_1\tprod\alpha_2\imply C}{
      \infer[\LR{\imply}]{\alpha_1\tprod\alpha_2, A, B, \G_2 \drv C}{
      \alpha_1, \alpha_2, A, B, \G_2 \drv C
      }
    }
  }
  &
  \infer[\LR{\imply}]{\alpha_1, \alpha_2, \alpha_1\tprod\alpha_2\imply C}{
    \infer[\RR{\tprod}]{\alpha_1, \alpha_2 \drv \alpha_1\tprod\alpha_2}{
      \alpha_1 \drv \alpha_1
      &
      \alpha_2 \drv \alpha_2
    }
    &
    C\drv C
  }
}
$$

or

{\small
$$
\hspace*{-3em}
\infer[\Cut]{\alpha_1, \alpha_2, A\tprod B, \G_2 \drv C}{
  \infer[\LR{\tprod}]{A\tprod B, \G_2 \drv \alpha_1\tprod\alpha_2\imply C}{
    \infer[\RR{\imply}]{A, B, \G_2 \drv \alpha_1\tprod\alpha_2\imply C}{
      \infer[\LR{\imply}]{\alpha_1\tprod\alpha_2, A, B, \G_2 \drv C}{
        \alpha_1, \alpha_2, A, B, \G_2 \drv C
      }
    }
  }
  &
  \infer[\Cut]{\alpha_1, \alpha_2, \alpha_2\imply\alpha_1\imply C\drv C}{
    \infer[\LR{\imply}]{\alpha_2, \alpha_1\imply\alpha_1\imply C\drv \alpha_1\imply C}{
      \alpha_2\drv\alpha_2
      &
      \infer[\RR{\imply}]{\alpha_1\imply C\drv\alpha_1\imply C}{
        \infer[\LR{\imply}]{\alpha_1, \alpha_1\imply C\drv C}{
          \alpha_1\drv \alpha_1
          &
          C\drv C
        }
      }
    }
    &
    \infer[\LR{\tprod}]{\alpha_1, \alpha_1\imply C\drv C}{
      \alpha_1\drv \alpha_1
      &
      C\drv C
    }
  }
}
$$
}
\item[Case 2]%
The bottom inference rule in $\Pi$ is$\LR{\dsum}$ \\[0.5em]
i.e.\;\;
$
\infer[\LR{\dsum}]{\alpha_1, \alpha_2, A \dsum B, \G_2 \drv C}{
  \alpha_1, \alpha_2, A, \G_2 \drv C & \G_1, B, \G_2 \drv C}
$ \\[1em]
Then we can construct a proof of \; 
$\alpha_1, \alpha_2, A \dsum B, \G_{2} \drv C$ in $\FL{}{}'$ as follows:
{\footnotesize
$$
\hspace*{-3em}
\infer[\Cut]{\alpha_1, \alpha_2, A\tprod B, \G_2 \drv C}{
  \infer[\LR{\dsum}]{A\dsum B, \G_2 \drv \alpha_1\tprod\alpha_2\imply C}{
    \infer[\RR{\imply}]{A, \G_2 \drv \alpha_1\tprod\alpha_2\imply C}{
      \infer[\LR{\tprod}]{\alpha_1\tprod\alpha_2, A, \G_2 \drv C}{
        \alpha_1, \alpha_2, A, \G_2 \drv C
      }
    }
    &
    \infer[\RR{\imply}]{B, \G_2 \drv \alpha_1\tprod\alpha_2\imply C}{
      \infer[\LR{\tprod}]{\alpha_1\tprod\alpha_2, B, \G_2 \drv C}{
        \alpha_1, \alpha_2, B, \G_2 \drv C
      }
    }
  }
  &
  \infer[\LR{\imply}]{\alpha_1, \alpha_2, \alpha_1\tprod\alpha_2\imply C}{
    \infer[\RR{\tprod}]{\alpha_1, \alpha_2 \drv \alpha_1\tprod\alpha_2}{
      \alpha_1 \drv \alpha_1
      &
      \alpha_2 \drv \alpha_2
    }
    &
    C\drv C
  }
}
$$
}

or

{\tiny
$$
\hspace*{-5em}
\infer[\Cut]{\alpha_1, \alpha_2, A\dsum B, \G_2 \drv C}{
  \infer[\LR{\dsum}]{A\dsum B, \G_2 \drv \alpha_2\imply\alpha_1\imply C}{
    \infer[\RR{\imply}]{A, \G_2 \drv \alpha_2\imply\alpha_1\imply C}{
      \infer[\RR{\imply}]{\alpha_2, A, \G_2 \drv \alpha_1\imply C}{
        \alpha_1, \alpha_2, A, \G_2 \drv C
      }
    }
    &
    \infer[\RR{\imply}]{B, \G_2 \drv \alpha_2\imply\alpha_1\imply C}{
      \infer[\RR{\imply}]{\alpha_2, B, \G_2 \drv \alpha_1\imply C}{
        \alpha_1, \alpha_2, B, \G_2 \drv C
      }
    }
  }
  &
  \infer[\Cut]{\alpha_1, \alpha_2, \alpha_2\imply\alpha_1\imply C\drv C}{
    \infer[\LR{\imply}]{\alpha_2, \alpha_1\imply\alpha_1\imply C\drv \alpha_1\imply C}{
      \alpha_2\drv\alpha_2
      &
      \infer[\RR{\imply}]{\alpha_1\imply C\drv\alpha_1\imply C}{
        \infer[\LR{\imply}]{\alpha_1, \alpha_1\imply C\drv C}{
          \alpha_1\drv \alpha_1
          &
          C\drv C
        }
      }
    }
    &
    \infer[\LR{\tprod}]{\alpha_1, \alpha_1\imply C\drv C}{
      \alpha_1\drv \alpha_1
      &
      C\drv C
    }
  }
}
$$
}
\end{description}
The other inference rules are handled in a similar way.

The other direction is clear since each axiom and inference rule of our $\FL{}{}'$ is a particular case of those of  Ono's $\FL{}{}$.

\begin{flushright}
{\it q.e.d.}
\end{flushright}




%
%
%
%
%

\section{Cut, Parameter in inference rules, associativity and distributivity}


In this section, we present some examples of  proof-figures of \FL{}{}' which show 
how associative law and distributive law are involved with the cut-elimination process,
and show how difficult it is to eliminate the applications of the cut rule in these proof-figure.


 \ \\
 \ \\

\noindent
Associativity\;:\;
$
A\tprod(B\tprod C) \drv (A\tprod B)\tprod C
$ \\[2em]

\noindent
First of all, we present the following proof figure in $\FL{}{}'$, which contains an application of the cut rule.

{\small
$$
\infer[\LR{\tprod}]{A\tprod(B\tprod C) \drv (A\tprod B)\tprod C}{
  \infer[\Cut]{A, B\tprod C \drv (A\tprod B)\tprod C}{
    \infer[\LR{\tprod}]{B\tprod C \drv A \imply (A\tprod B)\tprod C}{
      \infer[\RR{\imply}]{B, C \drv A \imply (A\tprod B)\tprod C}{
        \infer[\RR{\tprod}]{A, B, C \drv (A\tprod B)\tprod C}{
          \infer[\RR{\tprod}]{A, B \drv A\tprod B}{A\drv A & B\drv B}
          &
          C\drv C
        }
      }
    } 
    &
    \infer[\LR{\imply}]{A, A \imply (A\tprod B)\tprod C \drv A \imply (A\tprod B)\tprod C}{
      A \drv A
      &
      A \imply (A\tprod B)\tprod C \drv A \imply (A\tprod B)\tprod C
    }
  }
}
$$
}

 \ \\

\noindent
This application of the cut rule becomes eliminable in $\FL{}{}$ as the next proof figure shows:

$$
\infer[\LR{\tprod}]{A\tprod(B\tprod C) \drv (A\tprod B)\tprod C}{
  \infer[\LR{\tprod}]{A, B\tprod C \drv (A\tprod B)\tprod C}{
    \infer[\RR{\tprod}]{A, B, C \drv (A\tprod B)\tprod C}{
          \infer[\RR{\tprod}]{A, B \drv A\tprod B}{A\drv A & B\drv B}
          &
          C\drv C
        }
  }
}
$$

 \ \\
 \ \\

\noindent

Associativity\;:\;
$
(A\tprod B)\tprod C \drv A\tprod(B\tprod C) 
$ \\[2em]

This direction of associativity can be proved in both $\FL{}{}'$ and $\FL{}{}$, 
as the next proof figure shows:

$$
\infer[\LR{\tprod}]{(A\tprod B)\tprod C \drv A\tprod(B\tprod C)}{
  \infer[\LR{\tprod}]{A\tprod B, C \drv A\tprod(B\tprod C)}{
    \infer[\LR{\tprod}]{A, B, C \drv A\tprod(B\tprod C)}{
      A \drv A
      &
      \infer[\RR{\tprod}]{B, C \drv B\tprod C}{B\drv B & C\drv C}
    }
  }
}
$$

 \ \\

\noindent
Distributivity\;:\;
$
A\tprod(B\dsum C) \drv (A\tprod B)\dsum(A\tprod C)
$ \\[2em]

In $\FL{}{}'$, we need an application of the cut rule to prove this direction of distributivity,
as the following proof figure shows:

{\tiny
$$
\hspace*{-3em}
\infer[\LR{\tprod}]{A\tprod(B\dsum C) \drv (A\tprod B)\dsum(A\tprod C)}{
  \infer[\Cut]{A, B\dsum C \drv (A\tprod B)\dsum(A\tprod C)}{
    \infer[\LR{\dsum}]{B\dsum C \drv A \imply((A\tprod B)\dsum(A\tprod C))}{
      \infer[\RR{\imply}]{B \drv A \imply((A\tprod B)\dsum(A\tprod C))}{
        \infer[\RR{\dsum}]{A,B\drv (A\tprod B)\dsum(A\tprod C)}{
          \infer[\RR{\tprod}]{A,B\drv A\tprod B}{A\drv A & B\drv B}
        }
      }
    &
      \infer[\RR{\imply}]{C \drv A \imply((A\tprod B)\dsum(A\tprod C))}{
        \infer[\RR{\dsum}]{A,C\drv (A\tprod B)\dsum(A\tprod C)}{
          \infer[\RR{\tprod}]{A,C\drv A\tprod C}{A\drv A & C\drv C}
        }
      }
    }
  &
    \infer[\RR{\imply}]{A, A\imply((A\tprod B)\dsum(A\tprod C))\drv (A\tprod B)\dsum(A\tprod C)}{
      A\drv A & (A\tprod B)\dsum(A\tprod C) \drv (A\tprod B)\dsum(A\tprod C)
    }
  }
}
$$
}
 \ \\

\noindent
This application of the cut rule becomes eliminable in $\FL{}{}$ as the next proof figure shows:

%

$$
\infer[\LR{\tprod}]{A\tprod(B\dsum C) \drv (A\tprod B)\dsum(A\tprod C)}{
  \infer[\RR{\dsum}]{A, B\dsum C \drv (A\tprod B)\dsum(A\tprod C)}{
    \infer[\RR{\dsum}]{A,B\drv (A\tprod B)\dsum(A\tprod C)}{
      \infer[\RR{\tprod}]{A,B\drv A\tprod B}{A\drv A & B\drv B}
    }
    &
    \infer[\RR{\dsum}]{A,C\drv (A\tprod B)\dsum(A\tprod C)}{
      \infer[\RR{\tprod}]{A,C\drv A\tprod C}{A\drv A & C\drv C}
    }
  }
}
$$

\noindent

Distributivity\;:\;
$
(A\tprod B)\dsum(A\tprod C) \drv A\tprod(B\dsum C) 
$ \\[2em]

This direction of distributivity is provable in both $\FL{}{}'$ and $\FL{}{}$.

$$
\infer[\LR{\dsum}]{(A\tprod B)\dsum(A\tprod C)\drv A\tprod(B\dsum C)}{
  \infer[\LR{\tprod}]{A\tprod B \drv A\tprod(B\dsum C)}{
    \infer[\RR{\tprod}]{A,  B \drv A\tprod(B\dsum C)}{
      A\drv A
      &
      \infer[\RR{\dsum}]{B\drv B\dsum C}{B\drv B}
    }
  }
  &
  \infer[\LR{\tprod}]{A\tprod C \drv A\tprod(B\dsum C)}{
    \infer[\RR{\tprod}]{A,  C \drv A\tprod(B\dsum C)}{
      A\drv A
      &
      \infer[\RR{\dsum}]{C\drv B\dsum C}{C\drv C}
    }
  }
}
$$

The above proofs indicate the way how associativity and distributivity matter for the non-eliminability of  applications of the cut rule in a given proof of non-commutative substructural logic.
Indeed, the above proof figures, showing how  associativity and distributivity are related to the cut rule, are obtained through the analysis of the (unsuccessful) reduction process for a non-normalizable proof in Gentzen style natural deduction for non-commutative substructural logic. In other words, the role of associativity and distributivity (in the process of "reduction") becomes clearer in the places where cut elimination fails.

It is an open problem whether cut elimination holds for \FL{}{}' 
if we add associativity and distributivity to \FL{}{}'.

\begin{small}

\end{small}

\begin{thebibliography}{10}
%
\bibitem{Nakatogawa1998}
K. Nakatogawa and T. Ueno.
\newblock On structural inference rules
  for Gentzen-style natural deduction, Part I.
\newblock In {\em Proceedings  of the Sixth Asian Logic Conference, 
              Beijing 1996}, pp. 199--221,World Scientific, 1998.

\bibitem{Ono1998}
H.~Ono.
\newblock Proof-theoretic methods in nonclassical logic --- an introduction.
\newblock In {\em Theories of Types and Proofs}, chapter~6. Mathematical
  Society of Japan, Tokyo, 1998.

\bibitem{Ueno1999}
T.~Ueno, O.~Watari, and K.~Nakatogawa.
\newblock On structural inference rules for {Gentzen}-style natural deduction,
  {Part II}(extended abstract).
\newblock In {\em The Seventh Asian Logic Conference Book of Abstracts},
  Hsi-Tou, Taiwan, June 1999.


\bibitem{Ueno2007}
T.~Ueno, O.~ Watari, and K.~Nakatogawa.
\newblock On Structural Inference Rules for {Gentzen}-style Natural Deduction,
 {Part II}.
\newblock In {\em Archive for Studies in Logic},
Vol. 8, no. 1, pp. 1-23, 2007. (\FL{}{} in \cite{Ueno2007} coresponds to \FL{}{}' in the present paper.)
URL:http://logic.let.hokudai.ac.jp/~koji/research/archive/UenoEtAl07Structural \\
InferenceRulesII.pdf



\bibitem{Wansing1993}
H.~Wansing.
\newblock {\em The Logic of Information Structures}.
\newblock Number 681 in Lecture Notes in Artificial Intelligence.
  Springer-Verlag, Berlin, 1993.

\bibitem{Ueno2000}
T.~Ueno.
\newblock Natural Deductions for Substructural Logics, PhD thesis.
\newblock Division of Mathematics, Hokkaido University, 2000, 3.


\bibitem{Watari2000}
O.~Watari, K.~Nakatogawa, and T.~Ueno.
\newblock Normalization theorems for substructural logics in {Gentzen}-style
  natural deduction, abstract of the talk at {2000 Annual Meeting of the
  Assocication for Symbolic Logic, University of Illinois at Urbana-Champaign,
  June 3--7, 2000}.
\newblock {\em The Bulletin of Symbolic Logic}, 6(3):390--391, Sep. 2000.
\end{thebibliography}
\end{document}